\documentclass[twocolumn]{autart}    
\usepackage{amsfonts}                
\usepackage{amssymb}

\usepackage[a4paper,margin={1.55cm,1.55cm,1.55cm,1.55cm}]{geometry}

\usepackage{graphicx}                

\usepackage{apalike}
\usepackage{color}

\usepackage{savesym}
\savesymbol{AND}
\usepackage{algorithm}
\usepackage{algorithmic}

\newtheorem{definition}{Definition}
\newtheorem{assumption}{Assumption}
\newtheorem{remark}{Remark}
\newtheorem{lemma}{Lemma}
\newtheorem{theorem}{Theorem}

\newcommand{\BEASN}{\begin{eqnarray*}}
\newcommand{\EEASN}{\end{eqnarray*}}
\newcommand{\BEAS}{\begin{eqnarray}}
\newcommand{\EEAS}{\end{eqnarray}}
\newcommand{\BEQ}{\begin{equation}}
\newcommand{\EEQ}{\end{equation}}
\newcommand{\BIT}{\begin{itemize}}
\newcommand{\EIT}{\end{itemize}}

\newcommand{\eg}{{e.g.}}
\newcommand{\ie}{{i.e.}}
\newcommand{\nn}{{\nonumber}}

\newcommand{\reals}{\mathbf{R}}
\newcommand{\expect}{\mathbb{E}}

\newcommand{\calX}{\mathcal{W}}

\newcommand{\sfA}{\mathbf{P}}
\newcommand{\bfg}{\mathbf{g}}
\newcommand{\bfgh}{\hat{\mathbf{g}}}
\newcommand{\bfx}{\mathbf{w}}
\newcommand{\bfy}{\mathbf{v}}
\newcommand{\bfz}{\mathbf{u}}
\newcommand{\bfr}{\mathbf{r}}
\newcommand{\hs}{\hspace{0em}}
\newcommand{\vc}{\textcolor[rgb]{0.00,0.00,0.00}}

\begin{document}
\setlength{\arraycolsep}{0.2em}

\begin{frontmatter}

\title{Optimal Distributed Stochastic Mirror Descent for Strongly Convex Optimization
\thanksref{footnoteinfo}}

\thanks[footnoteinfo]{This paper was not presented at any IFAC meeting.
This work was supported in part by the Natural Science Fund for Excellent Young Scholars of Jiangsu Province under Grant BK20170099, in part by the National Natural Science Foundation of China under Grants 61573344, 61733018 and 61374180, and in part by the Research Grants Council of the Hong Kong Special Administrative Region, China under Grant CityU 11300415.
Corresponding author D.~Yuan. Tel. +86 25 85866512.
Fax +86 25 85866512.}

\author[Nust,Paestum]{Deming Yuan}\ead{dmyuan1012@gmail.com},
\author[Rome]{Yiguang Hong}\ead{yghong@iss.ac.cn},
\author[Baiae]{Daniel W. C. Ho}\ead{madaniel@cityu.edu.hk},
\author[Paestum]{Guoping Jiang}\ead{jianggp@njupt.edu.cn},

\address[Nust]{School of Automation, Nanjing University of
Science and Technology, Nanjing 210094, Jiangsu, P.R. China}
\address[Paestum]{College of Automation, Nanjing University of
Posts and Telecommunications, Nanjing 210023, Jiangsu, P.R. China}
\address[Rome]{Key Laboratory of Systems and Control,
Academy of Mathematics and Systems Science,
Chinese Academy of Sciences, Beijing 100190, P.R. China}
\address[Baiae]{Department of Mathematics,
City University of Hong Kong, Kowloon, Hong Kong}


\begin{abstract}
In this paper we consider convergence rate problems for stochastic strongly-convex optimization in the non-Euclidean sense with a
constraint set over a time-varying multi-agent network.
We propose two efficient non-Euclidean stochastic subgradient descent algorithms based on the Bregman divergence as distance-measuring function
rather than the Euclidean distances that were employed by the standard distributed stochastic projected subgradient algorithms.
For distributed optimization of non-smooth and strongly convex functions
whose only stochastic subgradients are available,
the first algorithm recovers the best previous known rate of $O(\ln(T)/T)$
(where $T$ is the total number of iterations).
The second algorithm is an epoch variant of the first algorithm
that attains the optimal convergence rate of $O(1/T)$,
matching that of the best previously known centralized stochastic subgradient algorithm.
Finally, we report some simulation results to illustrate the proposed algorithms.
\end{abstract}

\end{frontmatter}

\section{Introduction}
Recent years have witnessed a growing interest in developing distributed
subgradient algorithms for solving convex constrained optimization problem,
where the objective function is the sum of the local convex objective functions of nodes in a network
(see, \eg,
Nedi\'{c}, Ozdaglar, \& Parrilo, 2010;
Zhu \& Mart\'{i}nez, 2012;
Lin, Ren, \& Song, 2016),
due to their widespread applications including sensor networks
(see, \eg,
Shi, Ling, Wu, \& Yin, 2015),
and smart grid
(see, \eg,
Yi, Hong, \& Liu, 2016;
Chang, Nedi\'{c}, \& Scaglione, 2014),
to name a few.

\vc{
Strong convexity has been widely studied in convex optimization, because strongly convex cost functions can be easily found in a variety of engineering application domains like sensor networks and
smart grids and strongly convex properties are actively used in regularization methods.
Take the ridge regression problem as an example, where the objective function
consists of the strongly convex Tikhonov regularization term for some performance improvement in optimization computation
(see, \eg, Shalev-Shwartz, \& Ben-David, 2014).
In light of the increasing attention to distributed optimization,
various distributed designs for optimizing strongly convex functions (in the Euclidean sense)
have been proposed in the literature (see Nedi\'{c} \& Olshevsky, 2016
and Tsianos \& Rabbat, 2012),
due to its wide application in many practical fields and
its potential to provide better guarantees of convergence performance.
}

Many algorithms have been developed over the past years to solve distributed convex optimization problems
(see, \eg,
Lu, Tang, Regier, \& Bow, 2011;
Liu, Qiu, \& Xie, 2014;
Chen \& Sayed, 2012;
Yuan, Ho, \& Xu, 2016;
Yuan, Ho, \& Hong, 2016;
Kia, Cort\'{e}s, \& Mart\'{i}nez, 2015;
Ram, Nedi\'{c}, \& Veeravalli, 2010).
Such algorithms require only the first-order information of
the objective functions and Euclidean projection onto the constraint set.
This makes the algorithms attractive for large-scale optimization problems.
\vc{
	Specifically, recently an $O(\ln T/\sqrt{T})$ rate of convergence has been established in
	Nedi\'{c} and Olshevsky (2015).
}
However, the aforementioned algorithms are inherently Euclidean,
in the sense that they rely on measuring distances based on Euclidean norms.
This means that it is challenging or infeasible to generate efficient projections for
certain objective functions and constraint sets,
taking the Euclidean projection onto the unit simplex as an example.
In this paper, we shall develop a class of distributed algorithms
that are built on mirror descent,
which generalizes the projection step using the Bregman divergence.
Bregman divergences are a general class of distance-measuring functions,
which include the Euclidean distance and Kullback-Leibler (KL) divergence as special cases.
\vc{
The work Xi, Wu, and Khan \cite{xi2014} presents a first study of the
distributed optimization algorithm that builds on mirror descent,
for solving the non-strongly and deterministic variant of problem (\ref{problem-ssco});
}
however, only convergence results are established for the proposed algorithm.

\vc{
Convergence rate is an important issue in the distributed design.
Although the aforementioned algorithms in the last paragraph can be applied to
distributed optimization of strongly convex functions,
it is desirable to develop algorithms by further exploiting the strongly convexity of the objective function, in order to provide better performance such as \emph{faster} convergence rates.
}
In Nedi\'{c} and Olshevsky \cite{nedic2016tac}, the authors proposed a distributed stochastic
subgradient-push algorithm for solving problem (\ref{problem-ssco}),
under the assumption that the stochastic gradients of the objective functions are Lipschitz.
In particular, the algorithm converges at an $O( \ln(T)/T )$ rate in the unconstrained case,
which is (to the best of our knowledge) the previously best known rate in the literature.
\vc{
The work Rabbat \cite{rabbat2015}
developed a distributed proximal subgradient algorithm,
that uses the Euclidean distance as the distance-measuring function,
for solving the unconstrained composite stochastic optimization problems;
they prove that the proposed algorithm converges at an $O( 1/T )$ rate,
under the smoothness assumptions on the objective functions.
}
The authors in Tsianos and Rabbat \cite{tsianos2012}
proposed a class of distributed algorithms (in both batch and online setting)
that converge at an $O( \ln(T)/T )$ rate in the constrained case,
without making the smoothness assumptions on the objective functions.
\vc{
	Notably, recently the work Lan, Lee, and Zhou (2017) proposed a class of distributed stochastic optimization algorithms that convergence at a rate of $O(1/T^2)$, however, note that the algorithms are built on the accelerated subgradient schemes that utilize two previous estimates in the subgraident step.
}

In this paper we focus on establishing the convergence rate of algorithms
for the distributed strongly convex constrained optimization problem in the following form
\BEQ
\begin{array}{lll}
\mathrm{minimize}     &  & F(\bfx) = \sum\limits_{i=1}^{m} F_{i}(\bfx)        \\
\mathrm{subject\ to}  &  & \bfx \in \calX
\end{array}
\label{problem-ssco}
\EEQ
where each $F_i$
is strongly convex in the non-Euclidean sense and maybe non-smooth,
and $\calX \subseteq \reals^d$ is a convex constraint set known to all the nodes in the network.
Moreover, the nodes can only compute the noisy subgradients
of their respective objective functions.
To be specific,
we assume that there exists a stochastic subgradient oracle,
which, for any point $\bfx\in\calX$, returns a random estimate $\bfgh_i(\bfx)$
of a subgradient $\bfg_i(\bfx) \in \partial F_i(\bfx)$ so that $\expect[\bfgh_i(\bfx)] = \bfg_i(\bfx)$,
where $\partial F_i(\bfx)$ denotes the subdifferential set of $F_i(\cdot)$ at $\bfx$.
It is well-known that for (centralized) stochastic optimization of
non-smooth and strongly convex functions, the optimal convergence rate is $O( 1/T )$
(see, \eg,
Hazan \& Kale, 2014).
This fact, combined with the above observations,
motivates us to consider the following questions:
1) Is it possible to develop a distributed stochastic mirror descent algorithm
that recovers the best previously known rate $O(\ln(T)/T)$,
for distributed optimization of non-smooth and strongly convex functions?
and 2) For the same optimization problem, is it possible to devise a variant of
the developed algorithm that attains the optimal $O(1/T)$ convergence rate?

In this paper, we give affirmative answers to the above questions.
Specifically, the main contributions of this paper are highlighted as follows:
\BIT
\item{}
We consider the construction of non-Euclidean algorithms
for distributed stochastic optimization of strongly convex functions
whose only stochastic subgradients are available.
The algorithms generalize the standard distributed
stochastic projected subgradient algorithms to the non-Euclidean setting.
Therefore, the proposed algorithms are more flexible,
in the sense that they enable us to generate efficient updates to better reflect the
geometry of the underlying optimization problem, by carefully choosing the Bregman divergence.
\item{}
We propose a distributed stochastic mirror descent (DSMD) algorithm to answer the first question.
In particular, we show that for a total number of $T$ iterations,
the proposed algorithm achieves an $O\left( \ln (T) / T \right)$ rate of convergence,
by exploiting the strongly convexity of the objective functions.
\vc{
The DSMD algorithm is a stochastic variant of the algorithm in Xi, Wu, and Khan \cite{xi2014}, where only asymptotic convergence is established.
}
In addition, this rate recovers the best previous known rate in
Nedi\'{c} and Olshevsky \cite{nedic2016tac}
and Tsianos and Rabbat \cite{tsianos2012}.
Moreover, in contrast to the algorithm in
Nedi\'{c} and Olshevsky \cite{nedic2016tac},
our proposed DSMD algorithm is in constrained setting,
which naturally arises in a number of applications where each node's estimate
has to lie within some decision space (see, \eg,
Nedi\'{c}, Ozdaglar, \& Parrilo, 2010).
\item{}
We propose an epoch variant of the DSMD algorithm, called Epoch-DSMD algorithm,
to answer the second question.
The Epoch-DSMD algorithm combines the strength of the epoch gradient descent algorithm
that is widely used in the machine learning community (see, \eg, Hazan \& Kale, 2014)
and the DSMD algorithm.
In particular, we prove by induction that the resulting point returned by the last epoch attains
the optimal $O( 1/T )$ rate of convergence,
which largely improve the $O\left( \ln (T) / T \right)$ rate obtained by Tsianos and Rabbat \cite{tsianos2012}
with the Euclidean norm. 
\EIT


\emph{Notation:}
Let $\reals^{d}$ be the $d$-dimensional vector space.
Write $\| \bfx \|_2$ to denote the Euclidean norm of a vector $\bfx\in\reals^d$, and $\left< \mathbf{w},\mathbf{v} \right>$ to denote the standard
inner product on $\reals^{d}$, for any $\mathbf{w},\mathbf{v}\in\reals^d$.
We denote by $[m]$ the set of integers $\{1,\dots,m\}$.
For a vector $\bfx$, we denote its $i$th component by $[\bfx]_i$.
We denote the $(i,j)$th element of a matrix $\sfA$ by $[\sfA]_{ij}$.
For a differentiable function $f$,
Let $\nabla f(\bfx)$ denote the gradient of $f(\cdot)$ at $\bfx$, and $\expect[X]$ denote the expected value of a random variable $X$.

\section{Problem Setting and Assumptions}
In this paper, we are interested in solving convergence rate problems for (\ref{problem-ssco})
over a time-varying multi-agent network.
Specifically, let $\mathcal{G}(t) = \left( \mathcal{V} , \mathcal{E}(t) , \mathbf{P}(t)\right)$
be a directed graph that represents the nodes' communication pattern at time $t$,
where $\mathcal{V} = \{ 1,\ldots,m \}$ is the node set,
$\mathcal{E}(t)$ is the set of activated links at time $t$,
and $\mathbf{P}(t)$ is the communication matrix at time $t$.
We make the following standard assumption on graph $\mathcal{G}(t)$
(see, \eg, Ram, Nedi\'{c}, \& Veeravalli, 2010;
Yuan, Ho, \& Xu, 2016).

\begin{assumption}\label{assump:matrix}
The graph $\mathcal{G}(t) = \left(\mathcal{V},\mathcal{E}(t), \mathbf{P}(t) \right)$
satisfies ($t=1,2,\ldots$):
\BIT
\item[(a)]
There exist a scalar $0<\xi<1$ such that $[\sfA(t)]_{ii}\geq\xi$ for all $i$ and $t$,
and $[\sfA(t)]_{ij}\geq\xi$ whenever $(j,i)\in\mathcal{E}(t)$;
\item[(b)]
$\sfA(t)$ is doubly stochastic, \ie, $\sum_{i=1}^{m} [\sfA(t)]_{ij} =1$
and $\sum_{j=1}^{m} [\sfA(t)]_{ij} =1$ for all $i$ and $j$;
\item[(c)]
There exists some positive integer $B$ such that the graph
$\left(\mathcal{V}, \bigcup_{t=sB+1}^{(s+1)B} \mathcal{E}(t) \right)$
is strongly connected for every $s\!\geq \!0$.
\EIT
\end{assumption}

We now give the definition of the Bregman divergence, which is crucial in developing the algorithms.
\begin{definition}
Given a strongly convex and differentiable distance-generating function
$\Phi: \reals^d \rightarrow \reals$, the \emph{Bregman divergence} induced by $\Phi$ is defined as follows:
\BEASN
D_{\Phi} (\bfx || \bfy) := \Phi(\bfx) - \Phi(\bfy)
- \left< \nabla \Phi(\bfy) , \bfx-\bfy \right>.
\label{def-Bregman}
\EEASN
\end{definition}


It is time to provide a non-Euclidean strong convexity assumption on problem (\ref{problem-ssco}).

\begin{assumption}\label{assump-bregman}
\BIT
\item[(a)]
$\Phi$ is $\sigma_\Phi$-strongly convex
with respect to the Euclidean norm, where (without loss of generality) $\sigma_\Phi \geq 1$,
\ie, for any two points $\bfx$, $\bfy\in\reals^d$,
\BEASN
\Phi (\bfx) &\geq& \Phi (\bfy) + \left< \nabla \Phi (\bfy) , \bfx -\bfy \right>
+ \frac{\sigma_\Phi}{2} \| \bfx -\bfy \|_2^2;
\EEASN
\item[(b)]
$F_i(\bfx)$ is $\sigma_F$-strongly convex with respect to function $\Phi$,
\ie, for any two points $\bfx$, $\bfy\in\calX$, we have for any $i=1,\ldots,m$,
\BEASN
F_i (\bfx) &\geq& F_i(\bfy) + \left< \bfg_{i}(\bfy) , \bfx -\bfy \right>
+ \sigma_F D_\Phi (\bfx || \bfy)
\EEASN
where \vc{$D_\Phi(\bfx||\bfy)$ are convex
in their second argument $\bfy$ for every fixed $\bfx$}.
\EIT
\end{assumption}

\begin{remark}
Note that Assumptions \ref{assump-bregman}(a) is standard in developing
mirror descent algorithms for solving convex optimization problems
(see, for example, Beck \& Teboulle, 2003; Xi, Wu, \& Khan, 2014).
Assumption \ref{assump-bregman}(b) is commonly used in developing distributed mirror descent algorithms for
strongly convex optimization (see, \eg, Xi, Wu, \& Khan, 2014).
In fact, under Assumption \ref{assump-bregman}(b), one can easily show that
each $F_i$ is also strongly convex with respect to the Euclidean norm,
by simply assuming the smoothness of function $\Phi$.
\end{remark}

The following assumption is about the stochastic subgradient, which has been widely used in
the literature (see, \eg, Nedi\'{c} \& Olshevsky, 2016).

\begin{assumption}\label{assump-bounded-grad}
At any point $\bfx \in \calX$, the stochastic subgradient of function $F_i$ satisfies:
\BEASN
\expect [\| \bfgh_{i}(\bfx) \|_2^2 ] \leq G^2, \qquad \forall i.
\EEASN
\end{assumption}
\begin{remark}
\vc{
The motivation of studying stochastic subgradients comes from the following considerations: i) in many situations only the noisy subgradients are available and it is easy for calculation; and ii) the stochastic subgradients can be used in reducing the cost by taking only one sample among multiple samples of the objective function of each agent.
}
\end{remark}

\section{Distributed Stochastic Mirror Descent}
\label{sec-sub-optimal}
In this section, we propose a distributed stochastic algorithm that is
built on mirror descent.
The details of the non-Euclidean algorithm are given in Algorithm 1.
Specifically, we will establish the explicit convergence rate of the proposed algorithm.
\begin{algorithm}
\caption{DSMD}
\label{alg-dsmd-NE}
\begin{algorithmic}[1]
\REQUIRE total number of iterations $T$ and step size sequence $\{\eta_{t}\}_{t=1}^{T}$
\ENSURE $\bfx_{i,1} \in \calX$ for all $i\in[m]$
\FOR{$t=1$ to $T$}
\STATE Query the stochastic subgradient oracle at $\bfx_{i,t}$
to get a random $\bfgh_{i,t} := \bfgh_{i}(\bfx_{i,t})$
\STATE Compute
\BEASN
\nabla \Phi (\bfy_{i,t}) &=& \nabla \Phi (\bfx_{i,t}) - \eta_t \bfgh_{i,t} \\
\bfz_{i,t+1} &=& \mathop{\arg\min}_{\bfx \in \calX}  D_{\Phi} (\bfx || \bfy_{i,t})  \\
\bfx_{i,t+1} &=& \sum_{j=1}^{m} [\sfA(t)]_{ij} \bfz_{j,t+1}
\EEASN
\ENDFOR
\end{algorithmic}
\end{algorithm}

Before presenting the main convergence results,
we provide two standard setups for the distributed stochastic mirror descent algorithms.
\BIT
\item[$\bullet$]
\emph{Euclidean setup:} $\Phi(\bfx) = \frac{1}{2} \| \bfx \|_2^2$, and the
associated Bregman divergence is $D_{\Phi}(\bfy || \bfx) = \frac{1}{2} \| \bfy - \bfx \|_2^2$.
\vc{
In this case the DSMD algorithm (i.e., Algorithm 1) reduces to a stochastic variant of the
distributed projected subgradient algorithm in Nedi\'{c}, Ozdaglar, and Parrilo \cite{nedic2010}.
}
\item[$\bullet$]
\emph{Simplex setup:} For this setup, suppose that our constraint set
$\calX = \Delta_d = \{ \bfx\in\reals^{d} : \sum_{i=1}^{d} [\bfx]_i = 1 , [\bfx]_i \geq 0, i\in[d] \}$.
Let $\Phi(\bfx) = \sum_{i=1}^{d} [\bfx]_i \ln  [\bfx]_i$,
and it is proved that $\Phi(\bfx)$ is $1$-strongly convex with respect to the $\ell_1$ norm
$\|\cdot\|_1$ over $\calX$ (see, \eg, Nemirovski, Juditsky, Lan, \& Shapiro, 2009).
Its associated Bregman divergence, also known as the Kullback-Leibler divergence, is
$
D_{\Phi}(\bfx || \bfy) = \sum_{i=1}^{d} [\bfx]_i \ln  \frac{[\bfx]_i}{[\bfy]_i}.
$
More importantly, in that case Step 3 in Algorithm 1
is equivalent to the following \emph{distributed stochastic entropic descent algorithm}:
\BEAS
\hs [\bfz_{i,t+1}]_j &\!=\!& \frac{ [\bfx_{i,t}]_j \exp(-\eta_t [\bfgh_{i,t}]_j ) }
{\sum_{\ell=1}^{d} [\bfx_{i,t}]_\ell \exp(-\eta_t [\bfgh_{i,t}]_\ell )} \quad j\in[d]
\label{alg-entropic-a} \\
\hs \bfx_{i,t+1} &\!=\!& \sum_{j=1}^{N} [\sfA(t)]_{ij} \bfz_{j,t+1}.
\label{alg-entropic-b}
\EEAS
Note that the projection step in the standard distributed projected subgradient algorithm
cannot be solved explicitly (in fact, it involves computing the solution of
$d$-dimensional nonlinear equation at each step, as pointed out by Beck and Teboulle \cite{beck2003orl}),
as opposed to the update (\ref{alg-entropic-a}).
\EIT

Now we establish the main result of the DSMD algorithm.
\begin{theorem}\label{theorem-NE}
Let Assumptions 1--3 hold.
Let $\eta_t = \frac{1}{\sigma_F t}$ for all $t=1,\ldots,T$
and $\bfx^{\ast} = \arg\min_{\bfx\in\calX} F(\bfx)$.
Then, for any $j\in[m]$ and $T\geq 3$, we have
\BEASN
\expect \left[\| \widehat{\bfx}_{j,T} - \bfx^{\ast} \|_2^2\right]
\leq \frac{c \ln (T)}{T} + \frac{c^{\prime}}{T}
\EEASN
where
$\widehat{\bfx}_{i,T} = \frac{1}{T} \sum_{t=1}^{T} \bfx_{i,t}$,
$c = \frac{2 G^2}{\sigma_F^2 \sigma_\Phi}
\left(1  + \frac{4 \alpha \beta m }{(1-\beta) \sigma_\Phi} \right)$
and
$c^{\prime} = \frac{2G}{\sigma_F \sigma_\Phi} \left(\frac{2\alpha\beta}{1-\beta} + 1 \right)
\sum_{i=1}^{m} \expect[ \|\bfx_{i,1}\|_2 ]$.
\end{theorem}

The proof of Theorem \ref{theorem-NE} relies on the following lemma, that establishes a bound on the differences
among the estimates of all the nodes in the network.

\begin{lemma}\label{lemma-NE-disagree}
Let Assumptions 1--3 hold.
Then
\BEASN
&&\sum_{t=1}^{T} \sum_{i=1}^{m} \expect[ \| \bfx_{i,t} - \bfx_{j,t} \|_2 ] \nn\\
&&\quad\leq m \left(\frac{2\alpha\beta}{1-\beta} + 1 \right)
\sum_{i=1}^{m} \expect[ \|\bfx_{i,1}\|_2 ]
+ \frac{2 \alpha \beta m^2 G}{(1-\beta) \sigma_\Phi } \sum_{t=1}^{T}  \eta_{t}.
\EEASN
\end{lemma}
\textbf{Proof.}
See Appendix \ref{app-lemma-NE-disagree}.
$\square $


Armed with Lemma \ref{lemma-NE-disagree},
we are ready to present the proof of Theorem \ref{theorem-NE}.

\textbf{Proof.}\
[Proof of Theorem \ref{theorem-NE}]
Let $\expect_{|t-1} [X]$ denote the expectation
conditioned on all the randomness until round $t-1$.
Hence, $\expect_{|t-1} [\bfgh_{i,t}] = \bfg_{i}(\bfx_{i,t})$.
This fact implies
\BEAS
\sum_{i=1}^{m} \expect_{|t-1} \left< \bfgh_{i,t} , \bfx_{i,t} \!-\! \bfx^{\ast} \right>
\!\!=\!\! \sum_{i=1}^{m}\! \left< \bfg_{i}(\bfx_{i,t}) ,
\bfx_{i,t} \!-\! \bfx^{\ast} \right>.
\label{theorem-NE-0}
\EEAS
Taking the total expectation and using  Assumption \ref{assump-bregman}(b), we have
\BEAS
&&\hs \sum_{i=1}^{m} \expect[ \left< \bfgh_{i,t} , \bfx_{i,t} - \bfx^{\ast} \right> ] \nn\\
&&\hs \quad\geq \sum_{i=1}^{m} \expect \left[  F_i(\bfx_{i,t}) - F_i(\bfx^{\ast})
+ \sigma_F D_\Phi (\bfx^{\ast} || \bfx_{i,t})  \right] \nn\\
&&\hs \quad\geq  \sum_{i=1}^{m} \expect [F_i (\bfx_{j,t})] - F(\bfx^{\ast})
+  \sigma_F \sum_{i=1}^{m}
\expect [ D_\Phi (\bfx^{\ast} || \bfx_{i,t}) ]  \nn\\
&&\hs \qquad- G \sum_{i=1}^{m}  \expect [ \| \bfx_{i,t} - \bfx_{j,t} \|_2 ]
\label{theorem-NE-1}
\EEAS
where the last inequality follows from the convexity of function $F_i$, that is,
$
F_i (\bfx_{i,t}) \geq F_i (\bfx_{j,t})
+ \left< \bfg_i(\bfx_{j,t})  , \bfx_{i,t} - \bfx_{j,t} \right>
\geq F_i (\bfx_{j,t}) - G  \| \bfx_{i,t} - \bfx_{j,t} \|_2
$,
because
$\| \bfg_i(\bfx_{j,t}) \|_2 = \| \expect [\bfgh_i(\bfx_{j,t})] \|_2
\leq  \expect [\| \bfgh_i(\bfx_{j,t}) \|_2 ] \leq
\left( \expect [\| \bfgh_i(\bfx_{j,t}) \|_2^2 ] \right)^{1/2} \leq G$.
On the other hand, by following an argument similar to that of Lemma 6 in Hazan and Kale (2014),
it is easy to show that
\BEAS
&&\sum_{i=1}^{m} \left< \bfgh_{i,t} , \bfx_{i,t} - \bfx^{\ast} \right>
\leq \frac{\eta_t}{2} \sum_{i=1}^{m}  \|  \bfgh_{i,t} \|_2^2 \nn\\
&&\qquad\qquad + \sum_{i=1}^{m} \frac{ D_\Phi (\bfx^{\ast} || \bfx_{i,t})
- D_\Phi (\bfx^{\ast} || \bfx_{i,t+1}) }{\eta_t}.
\label{theorem-NE-1b}
\EEAS
Combining inequalities (\ref{theorem-NE-1}) and (\ref{theorem-NE-1b}), we get
\BEAS
&&\sum_{t=1}^{T}  \expect [ F (\bfx_{j,t}) ] - F (\bfx^{\ast})  \nn\\
&&\quad\leq \sum_{t=1}^{T} \sum_{i=1}^{m}
\frac{\expect [D_\Phi (\bfx^{\ast} || \bfx_{i,t}) ]
- \expect[D_\Phi (\bfx^{\ast} || \bfx_{i,t+1}) ]}{\eta_t}  \nn\\
&&\qquad - \sigma_F \sum_{i=1}^{m} \expect [ D_\Phi (\bfx^{\ast} || \bfx_{i,t}) ]
+ \sum_{t=1}^{T} \frac{\eta_t}{2} \sum_{i=1}^{m}  \expect [\|  \bfgh_{i,t} \|_2^2 ] \nn\\
&&\qquad+ G \sum_{t=1}^{T}  \sum_{i=1}^{m}  \expect[ \| \bfx_{i,t} - \bfx_{j,t} \|_2 ] \nn\\
&&\quad:= p_1 + p_2 + p_3 + p_4
\label{theorem-NE-2}
\EEAS
where $p_1$, $p_2$, $p_3$, $p_4$ denote the respective right-hand side terms in (\ref{theorem-NE-2}).
For terms $p_1$ and $p_2$, due to $\eta_t = \frac{1}{\sigma_F t}$ and
the non-negativity of the Bregman divergence, we obtain
\BEAS
p_1 + p_2 &=& \left( \eta_1^{-1} - \sigma_F \right)
\sum_{i=1}^{m} \expect[ D_\Phi (\bfx^{\ast} || \bfx_{i,1}) ] \nn\\
&&+ \sum_{t=2}^{T} \left( \eta_t^{-1} - \eta_{t-1}^{-1} - \sigma_F \right)
\sum_{i=1}^{m} \expect[ D_\Phi (\bfx^{\ast} || \bfx_{i,t}) ] \nn\\
&&-  \frac{1}{\eta_T} \sum_{i=1}^{m} \expect[ D_\Phi (\bfx^{\ast} || \bfx_{i,T+1}) ]
\leq 0.
\label{theorem-NE-3}
\EEAS
Term $p_3$ can be bounded by,
$
\sum_{t=1}^{T} \frac{\eta_t}{2} \sum_{i=1}^{m}  \expect [ \|  \bfgh_{i,t} \|_2^2 ]
\leq \frac{mG^2}{2}  \sum_{t=1}^{T} \eta_t
$,
by using Assumption \ref{assump-bounded-grad}.
In addition, $p_4$ can be bounded by Lemma \ref{lemma-NE-disagree}.
Hence, combining Lemma \ref{lemma-NE-disagree} with the preceding estimates,
and using the inequalities that
$
\frac{1}{T} \sum_{t=1}^{T} \eta_t = \frac{1}{T} \sum_{t=1}^{T} \frac{1}{\sigma_F t}
\leq \frac{2}{\sigma_F} \cdot \frac{\ln (T)}{T}
$, $\forall T\geq 3$ and
\BEASN
&&\frac{1}{T} \sum_{t=1}^{T}  F\left( \bfx_{j,t}\right) - F(\bfx^{\ast}) \nn\\
&&\quad \geq \left< \bfg(\bfx^{\ast}) , \widehat{\bfx}_{j,T} - \bfx^{\ast} \right>
+ m \sigma_F D_\Phi (\widehat{\bfx}_{j,T} || \bfx^{\ast}) \nn\\
&&\quad \geq m \sigma_F D_\Phi (\widehat{\bfx}_{j,T} || \bfx^{\ast}) \nn\\
&&\quad \geq \frac{m \sigma_F \sigma_\Phi}{2} \| \widehat{\bfx}_{j,T} - \bfx^{\ast} \|_2^2
\EEASN
where the first inequality follows from the first order optimality condition, we complete the proof.
$\square $


\begin{remark}
\vc{
Note that our algorithm is motivated by the seminal work on distributed optimization Nedi\'{c} and Ozdaglar (2009) and Nedi\'{c}, Ozdaglar, and Parrilo (2010). In contrast to the work Nedi\'{c}, Ozdaglar, and Parrilo (2010) that built on the Euclidean projection, the DSMD algorithm is based on mirror descent, that generalizes the Euclidean projection step by using the Bregman distance. This means that the DSMD algorithm allows efficient projections by carefully choosing the Bregman divergence, taking the unit simplex constraint set as an example (see (\ref{alg-entropic-a})--(\ref{alg-entropic-b})). Moreover, we have established the explicit convergence rate for the DSMD algorithm, while in Nedi\'{c}, Ozdaglar, and Parrilo (2010) only asymptotic convergence is obtained. Note also that our DSMD algorithm is a stochastic variant of the algorithm in Xi, Wu, and Khan \cite{xi2014}, and we have established non-asymptotic convergence rate results for the proposed algorithm, while in Xi, Wu, and Khan \cite{xi2014} only asymptotic convergence is obtained.
}
\end{remark}



\section{The Epoch-based Optimal Algorithm}
\label{sec-optimal}
In the previous section,
we proposed a distributed stochastic mirror descent algorithm
to achieve a rate of convergence at $O(\ln (T) / T)$,
which is suboptimal for stochastic strongly-convex optimization.
In this section, we present an optimal distributed stochastic mirror descent algorithm,
called Epoch-DSMD, to solve problem (\ref{problem-ssco}) and analyze its convergence properties.

\begin{algorithm}
\caption{Epoch-DSMD}
\label{alg-dsmd-Epoch}
\begin{algorithmic}[1]
\REQUIRE an initial step size $\eta_1$, number of iterations in the first epoch $T_1$,
and total number of iterations $T$
\ENSURE $\bfx^{1}_{i,1} = \arg\min_{\bfx\in\calX} \Phi(\bfx)$ for all $i\in[m]$, and set $k=1$
\WHILE{$\sum_{\ell=1}^{k} T_\ell \leq T$}
\FOR{$t=1$ to $T_k$}
\STATE Query the stochastic subgradient oracle at $\bfx^{k}_{i,t}$
to get a random $\bfgh^{k}_{i,t} := \bfgh_{i}(\bfx^{k}_{i,t})$
\STATE Compute
\BEASN
\nabla \Phi (\bfy^{k}_{i,t}) &=& \nabla \Phi (\bfx^{k}_{i,t}) - \eta_k \bfgh^{k}_{i,t} \\
\bfz^{k}_{i,t+1} &=& \mathop{\arg\min}_{\bfx \in \calX}  D_{\Phi} (\bfx || \bfy^{k}_{i,t})  \\
\bfx^{k}_{i,t+1} &=& \sum_{j=1}^{m} [\sfA_{k}(t)]_{ij} \bfz^{k}_{j,t+1}
\EEASN
\ENDFOR
\STATE Compute $\widehat{\bfx}^k_i \!=\! \frac{1}{T_k}\! \sum_{t=1}^{T_k} \! \bfx^{k}_{i,t}$
and update $\bfx^{k+1}_{i,1} \!\!=\!\! \widehat{\bfx}^k_i$
\STATE Update $T_{k+1} = 2 T_{k}$ and $\eta_{k+1} = \frac{1}{2} \eta_{k}$
\STATE Update $k = k+1$
\ENDWHILE
\end{algorithmic}
\end{algorithm}


To present the main convergence results of the Epoch-DSMD algorithm,
we assume from now on that the constraint set has finite radius
$R_\calX = \max_{\bfx,\bfy\in\calX} \| \bfx - \bfy \|_2$,
and denote the $\Phi(\cdot)$-diameter of $\calX$ by
$R_{\Phi,\calX} := \left( \max_{\bfx\in\calX} \Phi(\bfx)
- \min_{\bfx\in\calX} \Phi(\bfx) \right)^{1/2}.$
\begin{theorem}\label{theorem-Epoch}
Under the conditions of Theorem \ref{theorem-NE},
set the parameters in Algorithm 2 as
$\eta_1 = \frac{1}{\sigma_F}$ and $T_1 = 4$.
Then the final estimates $\bfx^{k^{\dag} + 1}_{i,1}$ enjoys a convergence rate of
\BEASN
\sum_{i=1}^{m}
\expect[ \| \bfx^{k^{\dag} + 1}_{i,1} - \bfx^{\ast} \|_2^2 ]
&\leq& \frac{64 \widehat{c} }{T}
\EEASN
where $k^{\dag} = \left\lfloor \log_2 \left( \frac{T}{4} + 1 \right) \right\rfloor$
is the total number of epochs in Algorithm 2,
and $\widehat{c} = \max\left\{ \frac{\sigma_F c_1 + 4 c_2 }{4 \sigma_F^2 \sigma_\Phi} ,
\frac{m R^2_{\Phi,\calX}}{4 \sigma_\Phi} \right\}$
with $c_1 = m G \left(\frac{2\alpha\beta}{1-\beta} + 1 \right)
\left( m R_\calX + \sum_{i=1}^{m} \expect[\|\bfx^1_{i,1}\|_2] \right)$
and $c_2 = \frac{mG^2}{2}  + \frac{2 \alpha \beta m^2 G^2}{(1-\beta) \sigma_\Phi}$.
\end{theorem}




\textbf{Proof.}\
First, we derive the following basic convergence result for Algorithm 2:
\BEAS
&&\sigma_F \sum_{i=1}^{m}
\expect [ D_\Phi (\bfx^{\ast} || \bfx^{k+1}_{i,1}) ] \nn\\
&&\quad\leq  \frac{c_1 }{T_k} + c_2 \eta_k
+ \frac{1}{\eta_k T_k} \sum_{i=1}^{m}
\expect[ D_\Phi (\bfx^{\ast} || \bfx^k_{i,1})].
\label{theorem-Epoch-9}
\EEAS
By following an argument similar to that of Theorem \ref{theorem-NE},
we have the following bound:
\BEAS
&&\hs \sum_{i=1}^{m} \left< \bfgh^k_{i,t} , \bfx^k_{i,t} - \bfx^{\ast} \right>
\leq \frac{\eta_k}{2} \sum_{i=1}^{m}  \|  \bfgh^k_{i,t} \|_2^2 \nn\\
&&\hs \qquad\ +\frac{1}{\eta_k} \sum_{i=1}^{m}
\left( D_\Phi (\bfx^{\ast} || \bfx^k_{i,t})
- D_\Phi (\bfx^{\ast} || \bfx^k_{i,t+1})  \right).
\label{theorem-Epoch-2}
\EEAS
Since $\eta_k$ is keeping constant in each epoch $k$,
by the fact that the term
$\sum_{i=1}^{m} \left( D_\Phi (\bfx^{\ast} || \bfx^k_{i,t})
- D_\Phi (\bfx^{\ast} || \bfx^k_{i,t+1}) \right)$
we can get a telescopic sum when summing over $t=1$ to $t=T_k$ and obtain
\BEAS
&&\sum_{t=1}^{T_k} \sum_{i=1}^{m}
\expect[  \left< \bfgh^k_{i,t} , \bfx^k_{i,t} - \bfx^{\ast} \right> ] \nn\\
&&\quad = \sum_{t=1}^{T_k} \sum_{i=1}^{m}
\expect[  \left< \bfg_{i}(\bfx^{k}_{i,t}) ,
\bfx^k_{i,t} - \bfx^{\ast} \right> ] \nn\\
&&\quad\leq \frac{\eta_k}{2} \sum_{t=1}^{T_k}
\sum_{i=1}^{m}  \expect[ \|  \bfgh^k_{i,t} \|_2^2 ]
+\frac{1}{\eta_k} \sum_{i=1}^{m} \expect [ D_\Phi (\bfx^{\ast} || \bfx^k_{i,1}) ] \nn\\
&&\quad\leq \frac{m G^2}{2} \eta_k T_k
+ \frac{1}{\eta_k} \sum_{i=1}^{m}  \expect[ D_\Phi (\bfx^{\ast} || \bfx^k_{i,1}) ]
\label{theorem-Epoch-3}
\EEAS
where the equality follows from the same reasoning as that of (\ref{theorem-NE-0}).
Let us turn our attention to the left-hand side of the preceding inequality.
It follows from Assumption \ref{assump-bregman}(b) that
\BEASN
&&\hs\sum_{t=1}^{T_k} \sum_{i=1}^{m}
\expect[  \left< \bfg_{i}(\bfx^{k}_{i,t}) ,
\bfx^k_{i,t} - \bfx^{\ast} \right> ] \nn\\
&&\hs\quad\geq \sum_{t=1}^{T_k} \sum_{i=1}^{m}
\expect\left[ F_i(\bfx^k_{i,t}) - F_i(\bfx^{\ast})
+ \sigma_F D_\Phi (\bfx^{\ast} || \bfx_{i,t})\right] \nn\\
&&\hs\quad\geq
\sigma_F\!\! \sum_{t=1}^{T_k}\! \sum_{i=1}^{m}\!
\expect [ D_\Phi (\bfx^{\ast} || \bfx^k_{i,t})]
\!\!-\! \!G  \sum_{t=1}^{T_k} \!\sum_{i=1}^{m}\!
\expect[ \| \bfx^k_{i,t} \!\!-\!\! \bfx^k_{j,t} \|_2 ]
\label{theorem-Epoch-4a}
\EEASN
where the last inequality is based on the same reasoning as that of (\ref{theorem-NE-1})
and the fact that $\bfx^{\ast}$ is the minimizer of problem (\ref{problem-ssco}).
Combining the preceding two inequalities,
and diving both sides by $T_k$, we get
\BEAS
&&\frac{\sigma_F}{T_k} \sum_{t=1}^{T_k} \sum_{i=1}^{m}
\expect [ D_\Phi (\bfx^{\ast} || \bfx^k_{i,t})] \nn\\
&&\quad \leq \frac{m G^2}{2} \eta_k
+ \frac{1}{\eta_k T_k} \sum_{i=1}^{m}  D_\Phi (\bfx^{\ast} || \bfx^k_{i,1}) \nn\\
&&\qquad+ \frac{G}{T_k} \sum_{t=1}^{T_k} \sum_{i=1}^{m}
\expect[ \| \bfx^k_{i,t} - \bfx^k_{j,t} \|_2 ] .
\label{theorem-Epoch-5}
\EEAS
The last term on the right-hand side can be easily bounded by using Lemma \ref{lemma-NE-disagree}, that is,
\BEAS
&&\hs\sum_{t=1}^{T_k} \sum_{i=1}^{m} \expect[ \| \bfx^k_{i,t} - \bfx^k_{j,t} \|_2 ] \nn\\
&&\hs\ \leq\!  m \!\! \left(\frac{2\alpha\beta}{1-\beta} \!+\! 1 \right)
\!\! \left( m R_\calX \!+\! \widehat{\delta}_1(\bfx) \right)
\!+\! \frac{2 \alpha \beta m^2 G}{(1-\beta) \sigma_\Phi} \eta_{k}  T_k
\label{theorem-Epoch-6}
\EEAS
where $\widehat{\delta}_1(\bfx) = \sum_{i=1}^{m} \expect[\|\bfx^1_{i,1}\|_2]$ and
in the last inequality we used the compactness assumption of the set $\calX$.
Substituting the bound (\ref{theorem-Epoch-6}) into inequality (\ref{theorem-Epoch-5}) gives
\BEAS
&&\frac{\sigma_F}{T_k} \sum_{t=1}^{T_k} \sum_{i=1}^{m}
\expect [ D_\Phi (\bfx^{\ast} || \bfx^k_{i,t})] \nn\\
&&\quad= \frac{c_1}{T_k} + c_2 \eta_k
+ \frac{1}{\eta_k T_k} \sum_{i=1}^{m}
\expect[ D_\Phi (\bfx^{\ast} || \bfx^k_{i,1})] .
\label{theorem-Epoch-7}
\EEAS
Using the assumption on the convexity of the Bregman divergence $D_\Phi$,
the left-hand side can be further lower bounded by the following:
\BEASN
\frac{\sigma_F}{T_k} \sum_{t=1}^{T_k} \sum_{i=1}^{m}
\expect [ D_\Phi (\bfx^{\ast} || \bfx^k_{i,t})]
&\geq& \sigma_F \sum_{i=1}^{m}
\expect [ D_\Phi (\bfx^{\ast} || \widehat{\bfx}^k_{i}) ] \nn\\
&=& \sigma_F \sum_{i=1}^{m}
\expect [ D_\Phi (\bfx^{\ast} || \bfx^{k+1}_{i,1}) ]
\label{theorem-Epoch-8}
\EEASN
where the last equality follows from Step 6 in Algorithm 2.
This, combined with (\ref{theorem-Epoch-7}), yields the bound (\ref{theorem-Epoch-9}).	
Now, we show by induction that
\BEAS
&&\hs\sum_{i=1}^{m} \expect[ D_\Phi (\bfx^{\ast} || \bfx^k_{i,1})] \nn\\
&&\hs\quad \leq \max\! \left\{ \frac{1}{\sigma_F} \! \left( \frac{c_1}{T_1} \!+\! c_2 \eta_1 \right) ,
\frac{1}{4} m R^2_{\Phi,\calX} \right\} \cdot 2^{-(k-3)}.
\label{theorem-Epoch-10}
\EEAS
We first prove that it is true at $k=1$.
It follows from the definition of the Bregman divergence that
$\sum_{i=1}^{m} D_\Phi (\bfx^{\ast} || \bfx^1_{i,1})
= \sum_{i=1}^{m} \big( \Phi(\bfx^{\ast}) - \Phi(\bfx^1_{i,1})
- \left< \nabla \Phi(\bfx^1_{i,1}) , \right.$
$\left. \bfx^{\ast}\!-\!\bfx^1_{i,1} \right> \big)$.
Utilizing the fact that $\bfx^{k}_{i,1} = \arg\min_{\bfx\in\calX} \Phi(\bfx)$
and applying the first order optimality condition for the term
$\left< \nabla \Phi(\bfx^1_{i,1}) , \bfx^{\ast}-\bfx^1_{i,1} \right>$,
we have
$
\left< \nabla \Phi(\bfx^1_{i,1}) , \bfx^{\ast}-\bfx^1_{i,1} \right>
\geq 0
$, which yields
\BEASN
\sum_{i=1}^{m}\! \expect[ D_\Phi (\bfx^{\ast} || \bfx^1_{i,1}) ]
&\!\leq\!& \sum_{i=1}^{m} \expect[ \Phi(\bfx^{\ast}) \!-\! \Phi(\bfx^1_{i,1}) ]
\!\leq\! m R^2_{\Phi,\calX}.
\label{theorem-Epoch-12}
\EEASN
Hence, the base of the induction holds. Assuming the bound (\ref{theorem-Epoch-10})
is true for $k$, we now claim that it holds for $k+1$ as well,
by combining the fact that $T_k = T_1 2^{k-1}$ and $\eta_k = \eta_1 2^{-(k-1)}$
with (\ref{theorem-Epoch-9}):
\BEAS
&&\hs \sum_{i=1}^{m}
\expect [ D_\Phi (\bfx^{\ast} || \bfx^{k+1}_{i,1}) ] \nn\\
&&\hs\quad\leq \max\left\{ \frac{1}{\sigma_F}\! \left( \frac{c_1}{T_1} \!+\! c_2 \eta_1 \right) ,
\frac{1}{4} m R^2_{\Phi,\calX} \right\} \cdot 2^{-(k-1)} \nn\\
&&\hs\qquad+ \frac{1}{4} \max\left\{ \frac{1}{\sigma_F}\! \left( \frac{c_1}{T_1} \!+\! c_2 \eta_1 \right) ,
\frac{1}{4} m R^2_{\Phi,\calX} \right\} \cdot 2^{-(k-3)} \nn\\
&&\hs\quad = \max \left\{ \frac{1}{\sigma_F}\! \left( \frac{c_1}{T_1} \!+\! c_2 \eta_1 \right) ,
\frac{1}{4} m R^2_{\Phi,\calX} \right\} \cdot 2^{-(k-2)}
\label{theorem-Epoch-13}
\EEAS
where the first inequality follows from
$\eta_1 = \frac{1}{\sigma_F}$ and $T_1 = 4$
and the induction hypothesis for $\sum_{i=1}^{m}
\expect[ D_\Phi (\bfx^{\ast} || \bfx^k_{i,1})]$.
This shows that the relation (\ref{theorem-Epoch-10}) holds for all $k\geq 1$.
Hence, combining the strongly convexity of function $\Phi$ with the relation (\ref{theorem-Epoch-10}),
it is easy to show that
\BEAS
&&\sum_{i=1}^{m}
\expect[ \| \bfx^{\ast} - \bfx^k_{i,1} \|_2^2 ] \nn\\
&&\quad\leq  \frac{2}{\sigma_\Phi}
\max\left\{ \frac{1}{\sigma_F} \left( \frac{c_1}{T_1} + c_2 \eta_1 \right) ,
\frac{1}{4} m R^2_{\Phi,\calX} \right\} \cdot 2^{-(k-3)} \nn\\
&&\quad=
\widehat{c} \cdot 2^{-(k-4)}.
\label{theorem-Epoch-14}
\EEAS
On the other hand, from the stopping criterion in Algorithm 2,
the number of epochs is given by the largest value of $k$ such that
$\sum_{i=1}^{k} T_k \leq T$, that is,
$
\sum_{i=1}^{k} T_1 2^{i-1} = T_1 (2^{k} - 1) \leq T
$, which implies that the final epoch is given by
$
k^{\dag}
= \left\lfloor \log_2 \left( \frac{T}{4} + 1 \right) \right\rfloor
$. Applying the bound (\ref{theorem-Epoch-14}) to $k^{\dag} + 1$ we get
\BEASN
\sum_{i=1}^{m}
\expect[ \| \bfx^{\ast} - \bfx^{k^{\dag} + 1}_{i,1} \|_2^2 ]
&\leq& \widehat{c} \cdot 2^{-(k^{\dag} + 1 - 4)} \nn\\
&\leq& \widehat{c} \cdot 2^{ - \left( \log_2 \left(\frac{T}{4} + 1\right) - 4 \right)}
\leq \frac{64 \widehat{c} }{T}
\label{theorem-Epoch-15}
\EEASN
where in the second inequality we used the relation
$k^{\dag} + 1 \geq \log_2 \left(\frac{T}{4} + 1\right)$.
The proof is complete.
$\square $

\begin{remark}
Theorem \ref{theorem-Epoch} shows that Algorithm 2 converges
at an $O(1/T)$ rate, matching that of the best previously known centralized
stochastic subgradient algorithm (see, for example, Hazan \& Kale, 2014).
To the best of our knowledge,
our proposed algorithm is the first distributed algorithm that utilizes the idea of epoch gradient descent
to develop a distributed stochastic mirror descent algorithm
(different from the distributed stochastic subgradient algorithm in Tsianos and Rabbat \cite{tsianos2012}).
Moreover, it achieves the optimal rate of convergence for distributed stochastic strongly
constrained convex optimization, without assuming smoothness of the objective functions,
in contrast to the $O(\ln(T)/T)$ rate of convergence in Tsianos and Rabbat \cite{tsianos2012}
and Nedi\'{c} and Olshevsky \cite{nedic2016tac} with Euclidean norms.
\end{remark}

\begin{remark}
\vc{
It is worth noting that the step size is keeping constant in each epoch,
which means that each node does not need to coordinate the step size
with those of its neighbors, but at the end of each epoch.
Specifically, $O(\ln (T))$ coordinations of the step size among all the nodes are needed.
This makes the Epoch-DSMD algorithm much easier to implement in a distributed setting,
as opposed to those algorithms that use (coordinated) diminishing step size
(see, \eg, Nedi\'{c} \& Olshevsky, 2016; Tsianos \& Rabbat, 2012; Xi, Wu, \& Khan, 2014).
}
\end{remark}

\section{Simulation Results}
In this section,
we consider the following standard distributed estimation problem over sensor networks
(see, \eg, Rabbat \& Nowak, 2004; Nedi\'{c} \& Olshevsky, 2016):
\BEQ
\min_{\bfx \in \calX}  \sum_{i=1}^{m} a_i \| \bfx - \mathbf{b}_i \|_2^2
\label{SimuRe}
\EEQ
where $a_i > 0 $ and $\mathbf{b}_i \in \reals^d$
are problem data known only to node $i$,
and $\calX$ is the constraint set known to all the nodes.
We will consider two constraint sets, namely,
(i) $\calX = \{ \bfx\in\reals^{d} : \sum_{i=1}^{d} [\bfx]_i = 1 , [\bfx]_i \geq 0, i\in[d] \}$;
and (ii) $\calX = \{ \bfx\in\reals^{d} : l_i \leq [\bfx]_i \leq u_i , i\in[d]\}$.

\vc{
Implement the proposed algorithms over a ring network that consists of $40$ nodes,
and the nodes are connected to form a single cycle. The network is time-varying,
in the sense that at each time instant, half of the links are activated randomly.
}
In all cases we use the dimension of estimate $d = 10$.
Note that the subgradient noises are random variables generated independent and
identically distributed from the normal distribution $\mathcal{N} (0, \sigma I_{d\times d})$,
where $\sigma$ is the magnitude of the noises that will be specified in the sequel.
\vc{
The simulation results for the methods are based on the average of 50 realizations.
}

\vc{
For the first constraint, \ie, the unit simplex,
the DSMD algorithm is just the distributed
stochastic entropic descent algorithm (\ref{alg-entropic-a})--(\ref{alg-entropic-b}).
We compare the convergence of our Epoch-DSMD algorithm with that of Rabbat \cite{rabbat2015}
(\ie, the MAMD algorithm).
Figs. 1 and 2 provide respectively a plot of the average error (on a log-scale)
versus the number of iterations $T$ for three randomly selected nodes,
for choices of $\sigma = 0.25$ and $\sigma = 0.5$.
It can be seen that both algorithms converge,
and Epoch-DSMD converges faster than MAMD.
Moreover, note that the MAMD algorithm involves a Euclidean projection onto the simplex at every step,
which, as we have stated earlier, is equivalent to computing the solution of a $d$-dimensional nonlinear equation.
This makes our proposed algorithms more favorable for this case.
}

\vc{
For the second constraint, \ie, a box constraint, we set the parameters in $\calX$ as follows:
$l_i = -1$, $u_i = 1$, for all $i$.
We compare the convergence of our Epoch-DSMD algorithm with that of
Nedi\'{c}, Ozdaglar, and Parrilo \cite{nedic2010} (\ie, the DSPS algorithm).
Note that in this case the DSMD algorithm reduces to a stochastic variant of the algorithm
in Nedi\'{c}, Ozdaglar, and Parrilo \cite{nedic2010}.
Figs. 3 and 4 provide respectively a plot of the average error (on a log-scale)
versus the number of iterations $T$ for three randomly selected nodes,
for choices of $\sigma = 0.25$ and $\sigma = 0.5$.
It can be seen from Figs. 3 and 4 that the proposed Epoch-DSMD algorithm
converges much faster than the standard distributed stochastic projected subgradient algorithm.
}

\begin{figure}[htb]
\begin{center}
\rotatebox{360}{\scalebox{0.25}[0.25]{\includegraphics{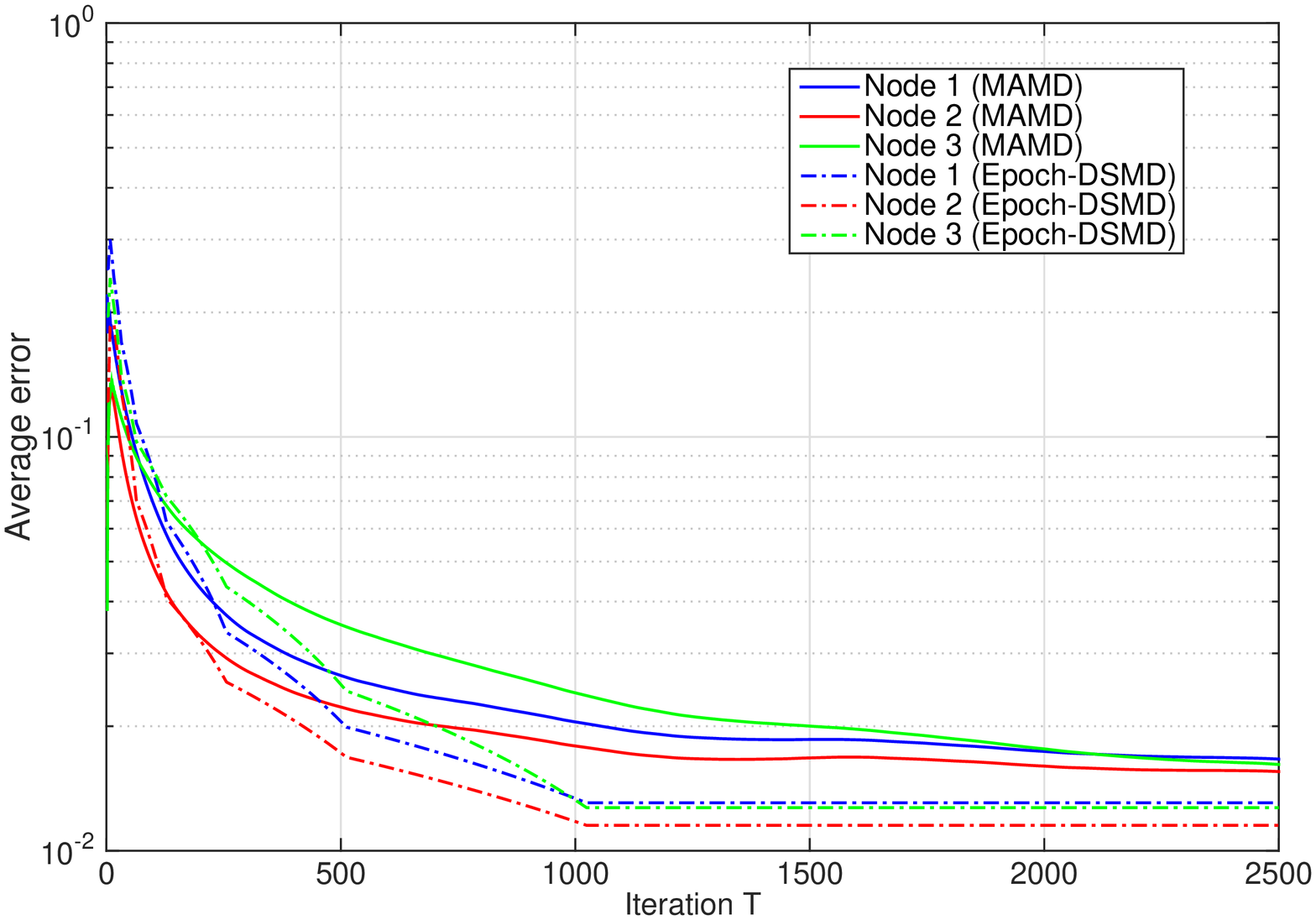}}}
\\[0pt]
{\normalsize {Fig. 1. The unit simplex constraint; $\sigma = 0.25$} }
\end{center}
\end{figure}

\begin{figure}[htb]
\begin{center}
\rotatebox{360}{\scalebox{0.25}[0.25]{\includegraphics{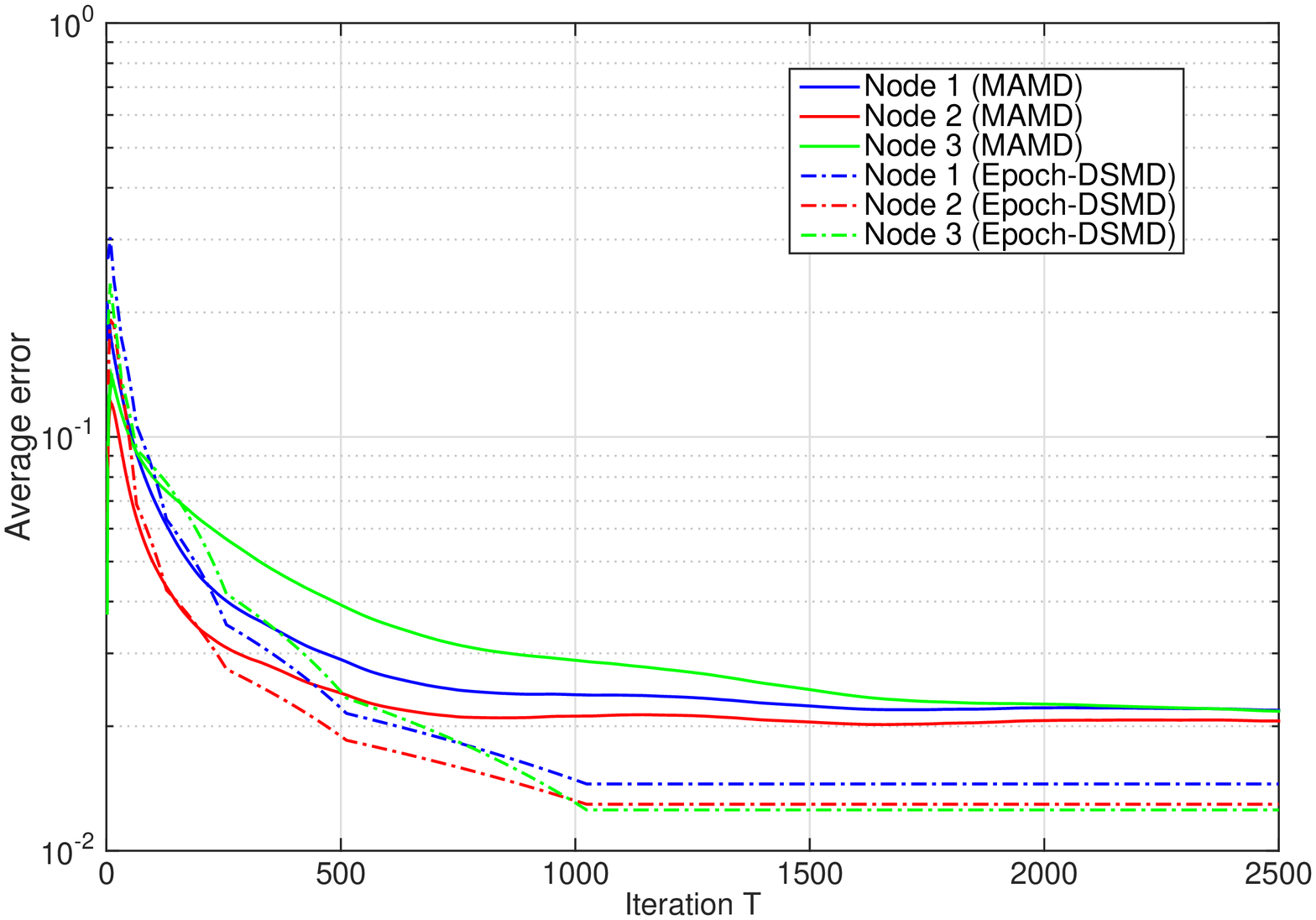}}}
\\[0pt]
{\normalsize {Fig. 2. The unit simplex constraint; $\sigma = 0.5$} }
\end{center}
\end{figure}

\begin{figure}[htb]
\begin{center}
\rotatebox{360}{\scalebox{0.25}[0.25]{\includegraphics{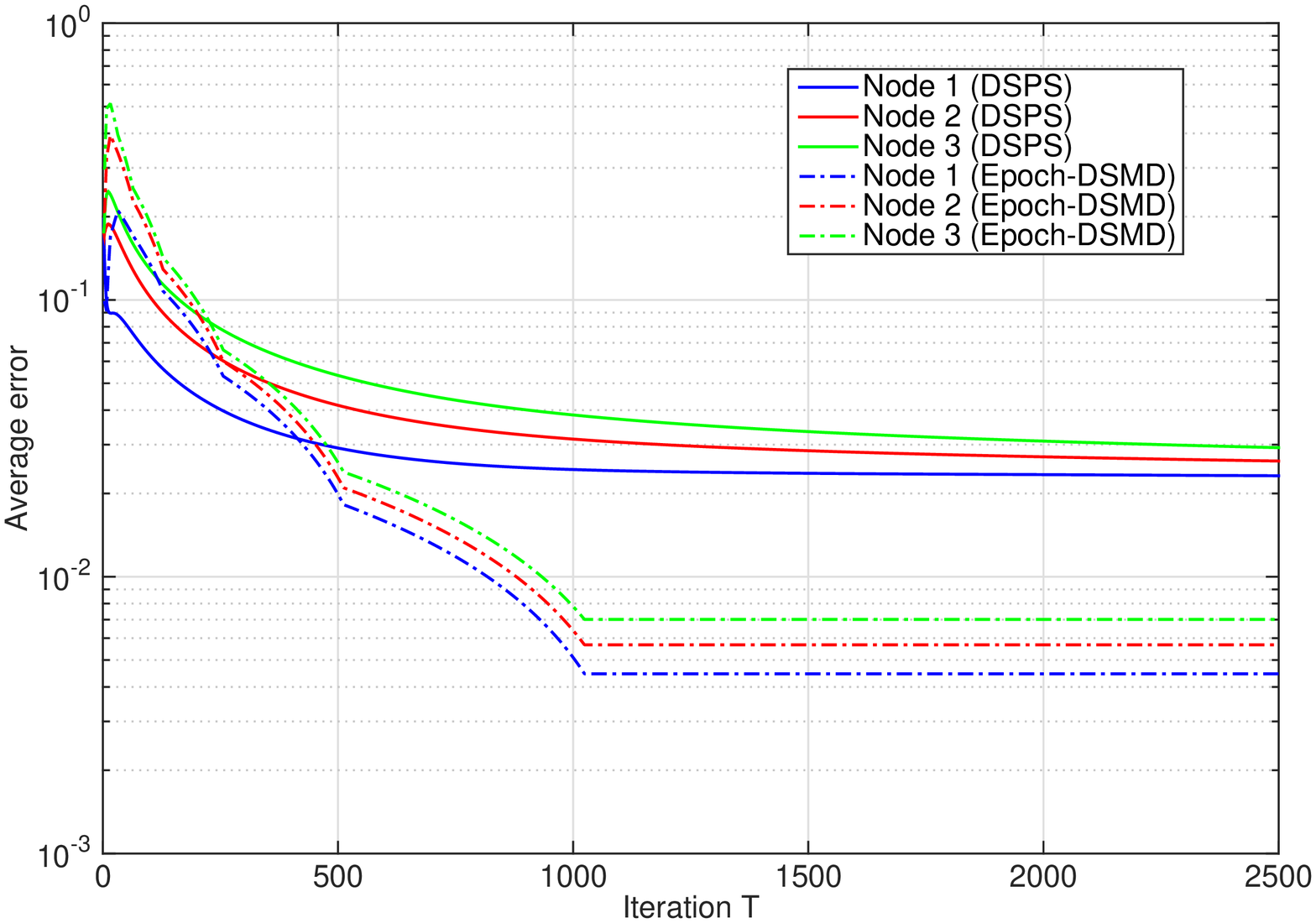}}}
\\[0pt]
{\normalsize {Fig. 3. The box constraint; $\sigma = 0.25$} }
\end{center}
\end{figure}

\begin{figure}[htb]
\begin{center}
\rotatebox{360}{\scalebox{0.25}[0.25]{\includegraphics{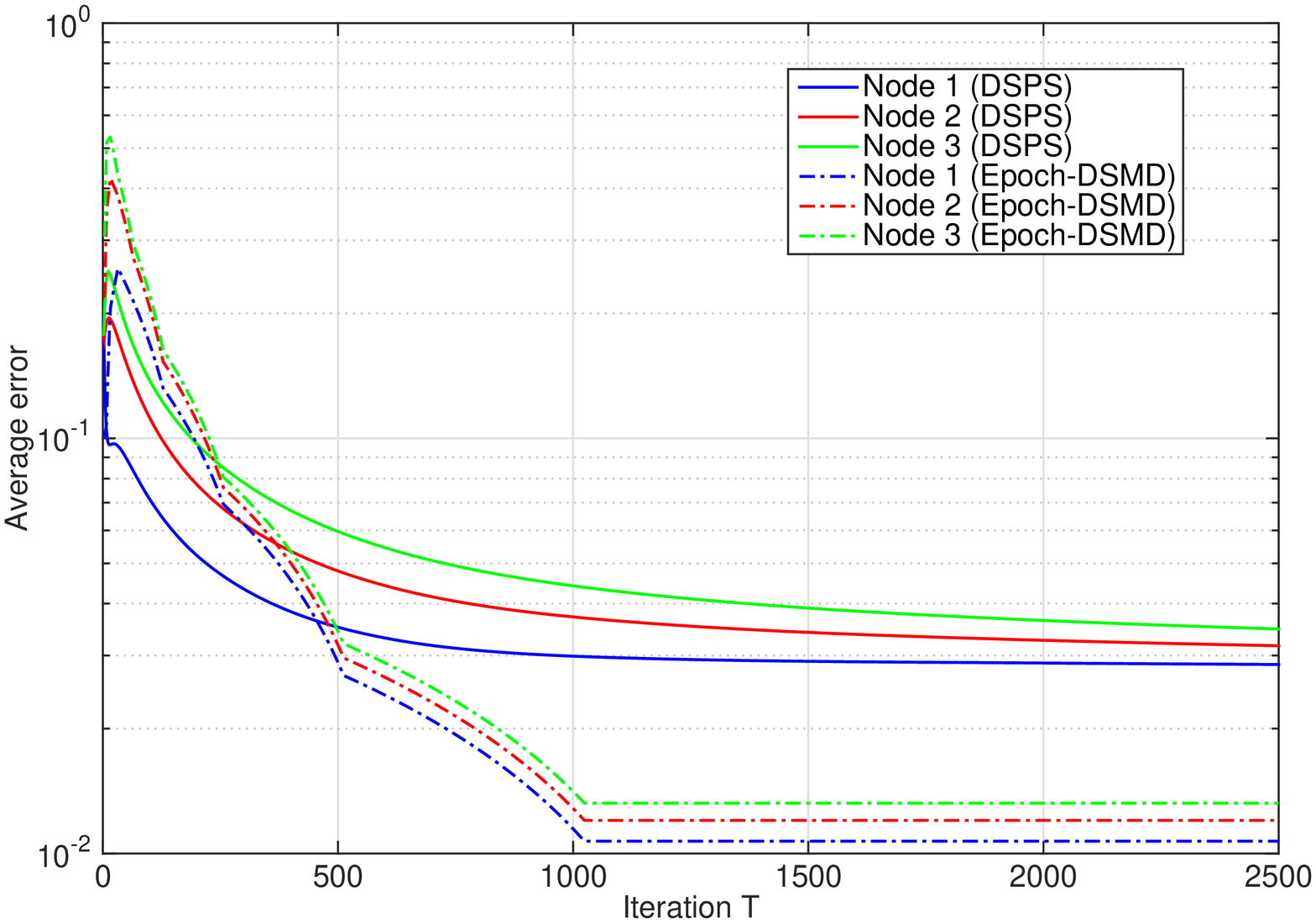}}}
\\[0pt]
{\normalsize {Fig. 4. The box constraint; $\sigma = 0.5$} }
\end{center}
\end{figure}

\section{Conclusion}
In this paper, we have studied the problem of distributed optimization of non-smooth and
strongly convex functions.
We have proposed two efficient non-Euclidean algorithms based on mirror descent.
The first algorithm recovers the best previously known rate,
and our second algorithm attains the optimal convergence rate.
\vc{
	There are several interesting questions that remain to be explored. For instance, one possible future research direction is to study the asynchronous variants of the proposed algorithms or remove the doubly stochasticity assumption on the weight matrix. Also, it would be of interest to adapt the accelerated subgradient schemes to the proposed algorithms to achieve an even faster convergence rate.
}

\appendix
\section{Proof of Lemma \ref{lemma-NE-disagree}}
\label{app-lemma-NE-disagree}
To simplify the notations, we denote
\BEAS
\overline{\bfx}_t = \frac{1}{m} \sum_{i=1}^{m} \bfx_{i,t}  \qquad\mathrm{and}\qquad
\bfr_{i,t} =    \bfz_{i,t+1} -  \bfx_{i,t}.
\label{lemma-NE-disagree-2}
\EEAS
Using the first order optimality condition for the update formula for $\bfz_{i,t+1}$
and noting that $\nabla D_\Phi (\bfz_{i,t+1} || \bfy_{i,t})
= \nabla \Phi (\bfz_{i,t+1}) - \nabla \Phi (\bfy_{i,t})$,
we obtain that, for all $\bfy\in\calX$,
$
\left< \nabla \Phi (\bfz_{i,t+1}) - \nabla \Phi (\bfy_{i,t}),
\bfz_{i,t+1} - \bfy\right> \leq 0.
$
Since $\bfx_{i,t} \in \calX$, by setting $\bfy = \bfx_{i,t}$ we have
\BEASN
\left< \nabla \Phi (\bfz_{i,t+1}) - \nabla \Phi (\bfy_{i,t}),
\bfz_{i,t+1} - \bfx_{i,t}\right> \leq 0.
\EEASN
Substituting the update formula for $\nabla \Phi  (\bfy_{i,t})$
into the preceding inequality yields
\BEASN
\eta_t\! \left< \bfgh_{i,t} , \bfx_{i,t} \!\!-\!\! \bfz_{i,t+1} \right>
&\!\geq\!&  \left< \nabla \Phi (\bfz_{i,t+1}) \!\!-\!\! \nabla \Phi (\bfx_{i,t}),
\bfz_{i,t+1} \!\!-\!\! \bfx_{i,t}\right>  \nn\\
&\!\geq\!& \sigma_{\Phi } \| \bfz_{i,t+1} \!- \! \bfx_{i,t} \|_2^2
\EEASN
because $\Phi $ is $\sigma_\Phi $-strongly convex.
This leads to the following bound
\BEAS
\expect [\| \bfr_{i,t} \|_2]
\leq \frac{\eta_t }{\sigma_\Phi } \expect[\| \bfgh_{i,t} \|_2]
\leq \frac{G}{\sigma_\Phi } \eta_t
\label{lemma-NE-disagree-3}
\EEAS
because $\expect [\| \bfgh_{i,t} \|_2 ] \leq
\left( \expect [\| \bfgh_{i,t} \|_2^2 ] \right)^{1/2} \leq G$,
based on Assumption \ref{assump-bounded-grad} and Jensen's inequality.

It is easy to derive the  general evolution of the average decision $\overline{\bfx}_t$,
by using the doubly stochasticity of the weight matrix $\sfA(t-1)$ and the definition of $\bfr_{i,t}$, that is
\BEASN
\overline{\bfx}_t = \overline{\bfx}_{1}
+  \sum_{\ell=1}^{t-1} \frac{1}{m} \sum_{i=1}^{m} \bfr_{i,\ell}  .
\label{lemma-disagree-bound-4}
\EEASN
Similarly, we derive the recursive relation for $\bfx_{i,t}$, and obtain
\BEASN
\bfx_{i,t} \!=\! \sum_{j=1}^{m} [\sfA(t-1,1)]_{ij} \bfx_{j,1}
+  \sum_{\ell=1}^{t-1}  \sum_{j=1}^{m} [\sfA(t-1,\ell)]_{ij} \bfr_{j,\ell}
\label{lemma-disagree-bound-5}
\EEASN
where $\sfA(t,\ell) = \sfA(t)\sfA(t-1)\cdots\sfA(\ell), \forall t\geq\ell\geq 1$.
Combining the preceding two equations,
and then taking the expectation and using the convergence property of the transition matrix $\sfA(t,\ell)$
\footnotemark[1]
\footnotetext[1]{
The convergence properties of $\sfA(t,\ell)$ can be characterized by the following lemma.\\
\textbf{Lemma}
\label{LEM:NEDIC}[Ram, Nedi\'{c}, and Veeravalli (2010)]
Under Assumption \ref{assump:matrix}, we have that, for all $i$, $j$ and all $t\geq \ell \geq 1$,
$\left| [\sfA(t,\ell)]_{ij} - 1/m \right| \leq \alpha \beta^{t-\ell+1}$,
with $\alpha = \left( 1 - \xi / 4m^2 \right)^{-2}$
and $\beta = \left( 1 - \xi / 4m^2 \right)^{1/B}$.
}
, we have that for any $t\geq 2$,
\BEASN
\expect[ \| \overline{\bfx}_t - \bfx_{i,t} \|_2 ]
&\leq&  \alpha \beta^{t-1} \delta_1(\bfx)
+  \sum_{\ell=1}^{t-1} \alpha \beta^{t-\ell} \sum_{j=1}^{m}
\expect [ \|\bfr_{j,\ell}\|_2 ] \nn\\
&\leq& \alpha \beta^{t-1} \delta_1(\bfx)
+   \frac{\alpha m G}{\sigma_\Phi}
\sum_{\ell=1}^{t-1}  \beta^{t-\ell}  \eta_\ell
\EEASN
where $\delta_1(\bfx) = \sum_{i=1}^{m} \expect [ \|\bfx_{i,1}\|_2 ]$,
and in the second inequality we used (\ref{lemma-NE-disagree-3}).
Summing the preceding inequalities over $t=1,\ldots,T$ and $i=1,\ldots,m$
and using some simple algebra,
we can derive the desired bound.


\begin{thebibliography}{99}


\bibitem[2009]{nedic2009a}
Nedi\'{c}, A., \& Ozdaglar, A. (2009).
Distributed subgradient methods for multi-agent optimization.
\emph{IEEE Transactions on Automatic Control},
\emph{54}(1), 48-61.


\bibitem[2010]{nedic2010}
Nedi\'{c}, A., Ozdaglar, A., \& Parrilo, P. A. (2010).
Constrained consensus and optimization in multi-agent networks.
\emph{IEEE Transactions on Automatic Control},
\emph{55}(4), 922--938.


\bibitem[2015]{cortes2015auto}
Kia, S. S., Cort\'{e}s, J., \& Mart\'{i}nez, S. (2015).
Distributed convex optimization via continuous-time coordination algorithms with discrete-time communication.
\emph{Automatica},
\emph{55}, 254--264.


\bibitem[2010]{ram2010}
Ram, S. S., Nedi\'{c}, A., \& Veeravalli, V. V. (2010).
Distributed stochastic subgradient projection algorithms for convex optimization.
\emph{Journal of Optimization Theory and Applications},
\emph{147}(3), 516--545.



\bibitem[2010]{cassandras2010tac}
Zhong, M., \& Cassandras, C. G. (2010).
Asynchronous distributed optimization with event-driven communication.
\emph{IEEE Transactions on Automatic Control},
\emph{55}(12), 2735--2750.


\bibitem[2012]{zhu2011tac}
Zhu, M., \& Mart\'{i}nez, S. (2012).
On distributed convex optimization under inequality and equality constraints.
\emph{IEEE Transactions on Automatic Control},
\emph{57}(1), 151--164.

\bibitem[2016]{lin2016auto}
Lin, P., Ren, W., \& Song, Y. (2016).
Distributed multi-agent optimization subject to nonidentical constraints and communication delays.
\emph{Automatica},
\emph{65}, 120--131.

\bibitem[2014]{liushuai2014}
Liu, S., Qiu, Z., \& Xie, L. (2014).
Continuous-time distributed convex optimization with set constraints.
In \emph{Proceedings of the 19th IFAC World Congress}
(pp. 9762--9767).


\bibitem[2011]{lu2011tac}
Lu, J., Tang, C. Y., Regier, P. R., \& Bow, T. D. (2011).
Gossip algorithms for convex consensus optimization over networks.
\emph{IEEE Transactions on Signal Processing},
\emph{56}(12), 2917--2923.

\bibitem[2014]{3s-book}
Shalev-Shwartz, S., \& Ben-David, S. (2014).
\emph{Understanding machine learning: From theory to algorithms},
New York, NY: Cambridge University Press.


\bibitem[2004]{rabbat2004}
Rabbat, M., \& Nowak, R. D. (2004).
Distributed optimization in sensor networks.
In \emph{Proceedings of the International Conference on Information Processing in Sensor Networks}
(pp. 20--27).



\bibitem[2016]{yuan2016cyber}
Yuan, D., Ho, D. W. C., \& Xu, S. (2016).
Regularized primal-dual subgradient method for distributed constrained optimization.
\emph{IEEE Transactions on Cybernetics},
\emph{46}(9), 2109--2118.

\bibitem[2016]{yuan2016siam}
Yuan, D., Ho, D. W. C., \& Hong, Y. (2016).
On Convergence rate of distributed stochastic gradient algorithm for convex optimization with inequality constraints.
\emph{SIAM Journal on Control and Optimization},
\emph{54}(5), 2872--2892.


\bibitem[2016]{hong2016auto}
Yi, P., Hong, Y., \& Liu, F. (2016).
Initialization-free distributed algorithms for optimal resource allocation with feasibility constraints and application to economic dispatch of power systems.
\emph{Automatica},
\emph{74}(2), 259--269.










\bibitem[2015]{shi2015siam}
Shi, W., Ling, Q., Wu, G., \& Yin, W. (2015).
EXTRA: an exact first-order algorithm for decentralized consensus optimization.
\emph{SIAM Journal on Optimization},
\emph{25}(2), 944--966.


\bibitem[2012]{sayed2012tsp}
Chen, J., \& Sayed, A. H. (2012).
Diffusion adaptation strategies for distributed optimization and learning over networks.
\emph{IEEE Transactions on Signal Processing},
\emph{60}(8), 4289--4305.


\bibitem[2014]{chang2014tac}
Chang, T. H., Nedi\'{c}, A., \& Scaglione, A. (2014).
Distributed constrained optimization by consensus-based primal-dual perturbation method.
\emph{IEEE Transactions on Automatic Control},
\emph{59}(6), 1524--1538.



\bibitem[2016]{nedic2016tac}
Nedi\'{c}, A., \& Olshevsky, A. (2016).
Stochastic gradient-push for strongly convex functions on time-varying directed graphs.
\emph{IEEE Transactions on Automatic Control},
\emph{61}(12), 3936--3947.


\bibitem[2015]{nedic2015tac}
\vc{
Nedi\'{c}, A., \& Olshevsky, A. (2015).
Distributed optimization over time-varying directed graphs.
\emph{IEEE Transactions on Automatic Control},
\emph{60}(3), 601--615.
}


\bibitem[2003]{beck2003orl}
Beck, A., \& Teboulle, M. (2003).
Mirror descent and nonlinear projected subgradient methods for convex optimization.
\emph{Operations Research Letters},
\emph{31}, 167--175.





\bibitem[2012]{tsianos2012}
Tsianos, K., \& Rabbat, M. (2012).
Distributed strongly convex optimization.
In \emph{Proceedings of the 50th Annual Allerton Conference on Communication, Control, and Computing}
(pp. 593--600).


\bibitem[2014]{xi2014}
Xi, C., Wu, Q., \& Khan, U. (2014).
Distributed mirror descent over directed graphs
(preprint, online at http://arxiv.org/abs/1412.5526).


\bibitem[2017]{lan2017}
\vc{
Lan, G., Lee, S., \& Zhou, Y. (2017).
Communication-efficient algorithms for decentralized and stochastic optimization
(preprint, online at https://www.arxiv.org/abs/1701.03961v2).
}



\bibitem[2015]{rabbat2015}
Rabbat, M. (2015).
Multi-agent mirror descent for decentralized stochastic optimization.
In \emph{Proceedings of the 6th International Workshop on Computational Advances in
Multi-Sensor Adaptive Processing (CAMSAP)} (pp. 517--520).

\bibitem[2014]{hazan2014-jmlr}
Hazan, E., \& Kale, S. (2014).
Beyond the regret minimization barrier: an optimal algorithm for stochastic strongly-convex optimization.
\emph{Journal of Machine Learning Research},
\emph{15}(1), 2489--2512.












\end{thebibliography}
\end{document}